\documentclass{amsart}
\usepackage{amsmath}
\usepackage{amsopn}
\usepackage{amssymb}
\usepackage{latexsym}
\usepackage{amsfonts}
\usepackage{amsthm}
\usepackage[all]{xy}

\newtheorem{thm}{Theorem}[section]
\newtheorem{pro}[thm]{Proposition}
\newtheorem{lem}[thm]{Lemma}
\newtheorem{cla}[thm]{Claim}

\newtheorem{cor}[thm]{Corollary}

\theoremstyle{definition}
\newtheorem{obs}[thm]{Observation}

\newtheorem{rem}[thm]{Remark}
\newtheorem{exa}[thm]{Example}
\newtheorem{defn}[thm]{Definition}

\newtheorem{conj}[thm]{Conjecture}

\newcommand{\een}{\end{enumerate}}
\newcommand{\blem}{\begin{lem}}
\newcommand{\elem}{\end{lem}}
\newcommand{\bcl}{\begin{cla}}
\newcommand{\ecl}{\end{cla}}
\newcommand{\ethm}{\end{thm}}
\newcommand{\bpr}{\begin{pro}}
\newcommand{\epr}{\end{pro}}
\newcommand{\bco}{\begin{cor}}
\newcommand{\eco}{\end{cor}}
\newcommand{\bcon}{\begin{conj}}
\newcommand{\econ}{\end{conj}}
\newcommand{\bde}{\begin{defn}}
\newcommand{\ede}{\end{defn}}
\newcommand{\bex}{\begin{exa}}
\newcommand{\eexa}{\end{exa}}
\newcommand{\bobs}{\begin{obs}}
\newcommand{\eobs}{\end{obs}}
\newcommand{\bexe}{\begin{exe}}
\newcommand{\eexe}{\end{exe}}

\begin{document}
\title[Congruence B-Orbits of Anti-Symmetric Matrices]{Involutions of the Symmetric Group and Congruence B-orbits of Anti-Symmetric Matrices}
\author[Y. Cherniavsky]{Yonah Cherniavsky}
\thanks{The author was supported by
the Swiss National Science Foundation.}
\address{Department of Mathematics and Computer Science, Ariel University Center of Samaria, Israel}
 \email{chrnvsk@gmail.com}

\begin{abstract} We present the poset of Borel congruence classes of anti-symmetric
matrices ordered by containment of closures. We show that there exists a  bijection between the set of these classes and the set of involutions of the symmetric group. We give two formulas for the rank function of this poset.
\end{abstract}
\maketitle
\section{Introduction}

The remarkable property of the Bruhat decomposition of $GL_n(\mathbb C)$ (i.e. the decomposition of $GL_n(\mathbb C)$ into the double cosets $\left\{B_1\pi B_2\right\}$ where $\pi\in S_n$ , $B_1,B_2\in\mathbb B_n(\mathbb C)$ -- the subgroup
of upper-triangular invertible matrices called the Borel subgroup) is that the natural order on double cosets (defined by the containment of closures) leads to the same poset as the combinatorially defined Bruhat order on permutations of $S_n$ (for $\pi,\sigma\in S_n$, $\pi\leqslant\sigma$ if $\pi$ is a subword of $\sigma$ with respect to the reduced form in Coxeter generators). L.~Renner introduced and developed the beautiful theory of Bruhat decomposition for not necessarily invertible matrices, see~\cite{R} and~\cite{R1}. When the Borel group acts on all the matrices, then the double cosets are in bijection with partial permutations which form a so called rook monoid $R_n$ which is the finite monoid whose elements are the 0-1 matrices with at most one nonzero entry in each row and column. The group of invertible elements of $R_n$ is isomorphic to the symmetric group $S_n$. Another efficient, combinatorial description of the Bruhat ordering on $R_n$ and a useful, combinatorial formula for the length function on $R_n$ are given by M.~Can and L.~Renner in~\cite{CR}.

The Bruhat poset of involutions of $S_n$ was first studied by F.~Incitti in~\cite{I} from purely combinatorial point of view. He proved that this poset is graded, calculated the rank function and also he showed several other important properties of this poset.
In~\cite{ABC} we give a geometric interpretation of the poset poset studied by F.~Incitti in~\cite{I} and its natural generalization considering the
action of the Borel subgroup on symmetric matrices by
congruence.

In this paper we present another graded poset of involutions of the symmetric group which also has the geometric nature. Denote by $\mathbb B_n(\mathbb C)$ the Borel
subgroup of $GL_n(\mathbb C)$, i.e. the group of invertible
upper-triangular $n\times n$ matrices over the complex numbers.
Denote by $\mathbb {AS}(n, \mathbb C)$ the set (which is actually a vector space with respect to standard operations of addition and multiplication by complex scalars, it is a Lie algebra usually denoted as $\mathfrak{so}$ with $[A,B]:=AB-BA$) of all complex anti-symmetric
$n\times n$ matrices. The congruence action of
$B\in\mathbb B_n(\mathbb C)$ on $S\in\mathbb S(n, \mathbb C)$ is
defined in the following way:
$S\,\,\longmapsto\,\,B^tSB\,\,.$
The orbits of this action (to be precisely correct we must say $S\,\mapsto\,\left(B^{-1}\right)^tSB^{-1}$ to get indeed a group action) are called  the {\it congruence B-orbits}. It is easy to see that $\mathbb {AS}(n, \mathbb C)$ is closed under this congruence action.

The main points of this paper are Proposition~\ref{asorball}, Definition~\ref{eqcount}, Theorem~\ref{PosetRankFunction} and Proposition~\ref{SecFm}. In Proposition~\ref{asorball} we show that the orbits of this action may be indexed by  involutions of $S_n$. In Definition~\ref{eqcount} we introduce the parameter $\mathfrak A$ and then in Theorem~\ref{PosetRankFunction} and Proposition~\ref{SecFm} we give two different formulas for the rank function of the studied poset using the parameter $\mathfrak A$. This parameter is similar to the parameter $\mathfrak D$ introduced in~\cite{ABC} and it can be seen as a particular case of a certain unified approach to the calculation of the rank function for several "Bruhat-like" posets as we briefly discuss it at the last section of~\cite{ABC}.

If we restrict this action on the set of invertible anti-symmetric matrices we get a poset of orbits that is isomorphic to the (reversed) Bruhat poset of involutions of $S_n$ without fixed points which is a subposet of the poset studied by F.~Incitti.
\section{A bijection between orbits and involutions}
The following Proposition~\ref{asorball} is somewhat similar to  Theorem~3.2 in~\cite{S}.
\bpr\label{asorball}
There is a bijection between the set of congruence B-orbits of all anti-symmetric $n\times n$ matrices and the set of all involutions of $S_{n}$.
\epr
\begin{proof} Let $A=\left(a_{ij}\right)_{i,j=1}^n\in\mathbb {AS}(n, \mathbb C)$. We move along the first row of $A$ until we find the first non-zero entry if there exists a non-zero entry in the first row, say $a_{1,j}\neq 0$. Notice that $j>1$ since all the diagonal entries of any anti-symmetric matrix are zeros. Now we can eliminate all non-zero entries in the first row from the right of $a_{1,j}$ multiplying $A$ from the right by appropriate upper-triangular matrix $B_1$ which differs from the identity matrix only in the $j$-th row:
$$
B_1=\begin{pmatrix}
1 &0 &\cdots &\cdots &\cdots &\cdots &\cdots &0\\
0 &1 &0 &\cdots &\cdots &\cdots &\cdots &0\\
\cdots &\cdots &\cdots &\cdots &\cdots &\cdots &\cdots &\cdots\\
0 &\cdots  &0 &1 &{-a_{1,j+1}\over a_{1,j}} &{-a_{1,j+2}\over a_{1,j}} &\cdots &{-a_{1,n}\over a_{1,j}}\\
\cdots &\cdots &\cdots &\cdots &\cdots &\cdots &\cdots &\cdots\\
0 &\cdots &\cdots &\cdots &0 &1 &0 &0\\
0 &\cdots &\cdots &\cdots &\cdots &0 &1 &0\\
0 &\cdots &\cdots &\cdots &\cdots &\cdots &0 &1
\end{pmatrix}\,.
$$

 Similarly we can eliminate all non-zero entries in the $j$-th column below $a_{1,j}$ multiplying $A$ from the left by appropriate lower-triangular matrix $C_1$:
 $$
C_1=\begin{pmatrix}
1 &0 &\cdots &\cdots &\cdots &\cdots &\cdots &\cdots &0\\
{-a_{2,j}\over a_{1,j}} &1 &0 &\cdots &\cdots &\cdots &\cdots &\cdots &0\\
\cdots &\cdots &\cdots &\cdots &\cdots &\cdots &\cdots &\cdots &\cdots\\
{-a_{j-1,j}\over a_{1,j}}&0 &\cdots &1 &0 &0 &\cdots &\cdots &0\\
0 &\cdots &\cdots &0 &1 &0 &0 &\cdots &0\\
{-a_{j+1,j}\over a_{1,j}}&0 &\cdots &0 &0 &1 &0 &\cdots &0\\
\cdots &\cdots &\cdots &\cdots &\cdots &\cdots &\cdots &\cdots &\cdots\\
{-a_{n-1,j}\over a_{1,j}} &\cdots &\cdots &\cdots &\cdots &\cdots &0 &1 &0\\
{-a_{n,j}\over a_{1,j}} &\cdots &\cdots &\cdots &\cdots &\cdots &\cdots &0 &1
\end{pmatrix}\,.
$$
Notice that there must be zero at the position $(j,1)$ in the matrix $C_1$ because there is zero at the position $(j,j)$ in the matrix $A$ since it is anti-symmetric. The matrix $C_1AB_1$ has zeros in the first row to the right from the position $(1,j)$ and in the $j$-th column below the position $(1,j)$. The matrix
$$C_1B_1^tAB_1C_1^t=\left(B_1C_1^t\right)^tAB_1C_1^t$$ is anti-symmetric, it is in the same congruence B-orbit as $A$ and its entries in the first row and $j$-th column are zeros except the entry $a_{1,j}$ in the position $(1,j)$, and also all the entries in the $j$-th row and first column are all zeros except the entry $-a_{1,j}$ in the position $(j,1)$.

Now we do the same elimination process for the second row and so on. At the end we get a monomial anti-symmetric matrix which is in the same congruence B-orbit as $A$. (Recall that a monomial matrix is a matrix which has at most one non-zero entry in each row and column.) Multiplying by the certain diagonal matrix from both sides we can get the monomial anti-symmetric matrix with only non-zero entries $\pm1$ and all minuses in the low triangle (below the main diagonal). Each such matrix encodes a the unique involution in the following way:  if we have $1$ at the position $(i,j)$, (and necessarily $-1$ at the position $(j,i)$), this pair of $\pm1$ corresponds to the transposition $(i,j)$ (the cycle of length 2) which exchanges $i$ and $j$. Since the matrix is anti-symmetric, if its $k$-th row is a zero row, then its $k$-th columns also is a zero column. So, if the $k$-th row and column are zeros we say that $k$ is a fixed point of the involution that we are constructing. Note that the zero matrix corresponds to the identity element of $S_{n}$ since everything is fixed.
\end{proof}
We illustrate the above proof with the following example:
\bex The monomial anti-symmetric matrix
 $\left[\begin{matrix}
0 &0 &0 &1 &0 &0\\
0 &0 &0 &0 &1 &0\\
0 &0 &0 &0 &0 &0\\
-1 &0 &0 &0 &0 &0\\
0 &-1 &0 &0 &0 &0\\
0 &0 &0 &0 &0 &0
\end{matrix}\right]$ corresponds to the involution\\ $\left(\begin{matrix}
1 &2 &3 &4 &5 &6\\
4 &5 &3 &1 &2 &6
\end{matrix}\right)\in S_6$, which can be written as the product of disjoint transpositions as $(1,4)(2,5)$.
\eexa
\bobs\label{asorbinv}
The congruence B-orbits of invertible anti-symmetric $2n\times 2n$ matrices can be indexed by involutions of $S_{2n}$ without fixed points.
\eobs
\noindent
\begin{proof} It is a particular case of Proposition~\ref{asorball}. We perform the elimination process described in the proof of Proposition~\ref{asorball}. Since the initial matrix is invertible, the monomial matrix $\pm1$'s which we get at the end is also invertible and therefore has no zero rows which means that the corresponding involution doesn't have fixed points.
\end{proof}
We end this section with an almost obvious observation:
\bobs
Let $X$ be an anti-symmetric matrix, $\pi$ a monomial anti-symmetric matrix and $B=\left(b_{ij}\right)$ invertible upper-triangular matrix such that $X=B^t\pi B$. Then $|Pf(X)|=|b_{11}b_{22}\cdots b_{nn}|$.
\eobs

\section{Partial order on orbits}

\label{OrbOrder}
When an algebraic group acts on a set of matrices, the classical
partial order on the set of all orbits is defined as follows:
$$\mathcal{O}_1 \leqslant_{\mathcal O}  \mathcal{O}_2\,\,\Longleftrightarrow\,\,\mathcal{O}_1 \subseteq \overline{\mathcal{O}_2}$$
where $\overline{S}$ is the  (Zarisski) closure
of the set $S$.

\subsection{Rank-control matrices}
\bde
 Let $X=\left(x_{ij}\right)$ be an
$n\times m$ matrix. For each $1 \leq k \leq n$ and $1 \leq l \leq m$, denote by $X_{k\ell}$ the
upper-left $k\times\ell$ submatrix of $X$. We denote by $R(X)$
the $n\times m$ matrix whose entries
are: $r_{k\ell}=rank\left(X_{k\ell}\right)$ and call it the {\it rank control matrix} of $X$.\ede
\begin{rem}
This rank-control matrix is an important tool in~\cite{ABC} and is similar to the one
introduced by A.~Melnikov~\cite{M} when she studied the poset (with
respect to the covering relation given in Definition~\ref{OrbOrder})
of adjoint B-orbits of certain nilpotent strictly upper-triangular
matrices.

The rank control matrix is connected also to the work of F.~Incitti~\cite{I} where the Bruhat poset of involutions of $S_n$ is
studied.
\end{rem}

 \bde\label{rankcontrolorder} Define the following order on
$n \times m$ matrices with positive integer entries: Let
$P=\left(p_{ij}\right)$ and $Q=\left(q_{ij}\right)$ be two such
matrices.

 Then
$$P\leqslant_{\mathcal R}
Q\,\,\Longleftrightarrow\,\,p_{ij}\leqslant q_{ij}\,\,\textrm{for
all}\,\,i,j\,.$$ \ede

The following proposition appears in another form as Theorem 2.1.5 in~\cite{BB}. Here we identify a permutation with the corresponding permutation matrix.

\bpr\label{BrOrdPerm} Denote by $\leqslant_B$
the Bruhat order of $S_n$ and let $\pi,\sigma\in S_n$. Then
$$\pi\leqslant_{\mathcal B}\sigma\quad\Longleftrightarrow\quad
R(\pi)\geqslant_{\mathcal R} R(\sigma)\,.$$ In other words, the
Bruhat order on permutation corresponds to the inverse order of
their rank-control matrices. \epr

\section{The Poset of Congruence B-Orbits of Anti-Symmetric Matrices}
The following easy proposition implies that the rank control matrix is an invariant of a congruence B-orbit.
\bpr\label{ulr} Let $X,Y\in GL_n(\mathbb F)$ be such that $Y=LXB$ for
some invertible lower-triangular matrix $L$ and some Borel (i.e.
invertible upper-triangular) matrix $B$. Denote by $X_{k\ell}$ and
$Y_{k\ell}$ the upper-left $k\times\ell$ submatrices of $X$ and $Y$
respectively. Then for all $1\leqslant k,\ell\leqslant n$
$$rank\left(X_{k\ell}\right)=rank\left(Y_{k\ell}\right)\,.$$
 \epr

Here is a direct consequence of Proposition~\ref{ulr}.
\bpr \label{BrCsRankCont}  All the matrices of a fixed
congruence B-Orbit have the same rank-control
matrix. In other words, if $X\in\mathbb{AS}(n,\mathbb C)$ and
$\mathcal A_X$ is the congruence B-orbit of $X$, then
$$
\mathcal A_X=\left\{S\in \mathbb {AS}(n,\mathbb
C)\,|\,R(S)=R(X)\right\}\,.
$$
\epr

Now we give the proposition which describes the orbit closures. This proposition follows from Theorem~15.31 given by E. Miller and B. Sturmfels, see~\cite[Chapter 15, page 301]{MS}:
\bpr\label{CosetClosure} Let $\pi$ be a partial involution and let
$R(\pi)$ be its rank-control matrix.  Then
$$
\overline{\mathcal A_X}=\left\{S\in \mathbb {AS}(n,\mathbb
C)\,|\,R(S)\leqslant_{\mathcal R} R(X)\right\}\,.
$$
 \epr
The next corollary characterizes the order relation of the poset of B-orbits.
\begin{cor}\label{main} Let $X,Y\in \mathbb {AS}(n,\mathbb
C)$.
Then
$$
\mathcal A_X\leqslant_{\mathcal O}\mathcal A_Y\iff R(X)
\leqslant_{\mathcal R}R(Y)
$$
\end{cor}

\section{The Rank Function of the Poset}
\bde A poset is called graded (or ranked) if for any two elements $x$ and $y$ of this poset any two maximal chains from $x$ to $y$ have the same length.
\ede
\bpr\label{rankdimrenner} The poset of congruence B-orbits of anti-symmetric matrices (with respect to the order $\leqslant_{\mathcal O}$) is a graded poset with the rank function given by the dimension of the closure.
\epr
\bco\label{rankdimrennerinv} The poset of congruence B-orbits of invertible anti-symmetric matrices (with respect to the order $\leqslant_{\mathcal O}$) is a graded poset with the rank function given by the dimension of the closure since it is an interval in the poset of congruence B-orbits of all anti-symmetric matrices.
\eco
\begin{rem}
It follows from Observation~\ref{asorbinv} and the results of~\cite{ABC} that the poset of congruence B-orbits of invertible anti-symmetric matrices which is the poset of involutions not having fixed points in the symmetric group is a graded subposet of the (reversed) Bruhat poset, while the whole poset of involutions of the symmetric group with the ordering given by the containment of closures of congruence B-orbits of anti-symmetric matrices is not the Bruhat poset.
\end{rem}
Proposition~\ref{rankdimrenner} is a particular case of the following fact. Let $G$ be a connected, solvable group acting on an irreducible, affine variety $X$. Suppose that there are a finite number of orbits. Let $O$ be the
set of $G$-orbits on $X$. For $x,y\in O$ define $x\leqslant y$ if $x\subseteq\overline y$. Then $O$ is a graded poset. This fact is given as an exercise in~\cite{R} (exercise 12, page 151) and can be proved using the proof of the Theorem of Section~8 of~\cite{R1}. (Our situation is a particular case of this fact because the Borel group is solvable, the variety of all anti-symmetric matrices is irreducible since it is a vector space and the number of orbits is finite since there are only finitely many involutions in the symmetric group.)

The question is to find an
algorithm which calculates the $\dim\overline{\mathcal A_X}$
from the involution of the symmetric group which corresponds to the orbit ${\mathcal A_X}$ or from the rank-control
matrix $R(X)$. Here we present such an algorithm.
\bde\label{eqcount} Let $X\in \mathbb {AS}(n,\mathbb
C)$ and let $R(X)=(r_{ij})_{i,j=1}^n$ be the rank-control matrix of $X$.  Add an extra
$0$ row to $R(X)$, pushed one place to the left, i.e. assume that
$r_{0k}=0$ for each $0\leqslant k < n$. Denote
$$
\mathfrak A(X)=\#\left\{(i,j)\,|\,1\leqslant i<
j\leqslant n\quad\textrm{and}\quad r_{ij}=r_{i-1,j-1}\right\}.
$$\ede
The parameter $\mathfrak A$ counts equalities in the diagonals of the upper triangle of the rank-control matrix and it is very similar to the parameter $\mathfrak D$ introduced in~\cite{ABC}. The only difference is that in the case of anti-symmetric matrices we consider the upper triangle without the main diagonal while in the case of the symmetric matrices which is studied in~\cite{ABC} the main diagonal is also considered. This difference is very natural since an anti-symmetric matrix is completely by its upper triangle without the main diagonal which consists of zeros while a symmetric matrix may have anything in its main diagonal.
From now on we don't distinguish between an involution $\pi$ of $S_n$ and the monomial anti-symmetric matrix whose non-zero entries are $\pm1$ (with minuses in the lower triangle) associated with $\pi$ by the bijection presented in Proposition~\ref{asorball}. Notice that when $\pi$ is a monomial matrix then the entry $r_{ij}$ of the rank-control matrix $R(\pi)$ is the number of nonzero entries of $\pi$ seen from the position $(i,j)$ when we are looking to the north-west.
\begin{thm}\label{PosetRankFunction} Let $\pi\in S_n$ be an involution. Then
$$
\dim\,\overline{\mathcal A_\pi}=\frac{n^2-n}{2}-\mathfrak A(\pi).
$$
\end{thm}
\begin{proof}
We have to explain the following before we begin the proof. By the variety which corresponds to some fragment of the matrix we mean the following: any variety of $n\times n$ matrices is a subset of the vector space $\mathbb C^{n^2}$ and the variety which corresponds to a fragment of $n\times n$ matrix is a projection the big variety on the corresponding subspace of $\mathbb C^{n^2}$. Denote by $V^{kn}$ the variety which corresponds to
$$
\left[\begin{matrix}
0 &a_{12} &\cdots &\cdots &a_{1,k} &\cdots &\cdots &a_{1,n-1} &a_{1,n} \\
-a_{12} &0 &\cdots&\cdots  &a_{2,k} &\cdots &\cdots &a_{2,n-1} &a_{2,n} \\
\cdots &\cdots &\cdots &\cdots&\cdots &\cdots &\cdots &\cdots &\cdots\\
-a_{1,k-1} &-a_{2,k-1} &\cdots &0 &a_{k-1,k}&\cdots &\cdots&a_{k-1,n-1} &a_{k-1,n}\\
-a_{1,k} &-a_{2,k} &\cdots &\cdots &0&\cdots &\cdots&a_{k,n-1}  &a_{k,n}\\
-a_{1,k+1} &-a_{2,k+1} &\cdots &\cdots&-a_{k,k+1} &0 &\cdots &a_{k+1,n-1} &\square\\
\cdots &\cdots &\cdots &\cdots &\cdots &\cdots&\cdots &\cdots&\square\\
-a_{1,n-1} &-a_{2,n-1} &\cdots &\cdots &\cdots &\cdots&\cdots &0 &\square\\
-a_{1,n} &-a_{2,n} &\cdots &\cdots&-a_{k,n}&\square&\square&\square&0
\end{matrix}\right]\,.
$$ (For $V^{kn}$ the last non empty entry in the $n$-th column is in the row number $k$). Consider also the variety $V^{k-1,n}$ which corresponds to
$$
\left[\begin{matrix}
0  &\cdots &a_{1,k-1}  &\cdots &a_{1,n-1} &a_{1,n} \\
\cdots  &\cdots&\cdots  &\cdots &\cdots &\cdots\\
-a_{1,k-1}  &\cdots &0 &\cdots&a_{k-1,n-1} &a_{k-1,n}\\
-a_{1,k}  &\cdots &-a_{k-1,k} &\cdots&a_{k,n-1}  &\square\\
\cdots  &\cdots &\cdots &\cdots &\cdots&\square\\
-a_{1,n-1}  &\cdots &\cdots  &\cdots&0 &\square\\
-a_{1,n}  &\cdots&-a_{k-1,n} &\square&\square&0
\end{matrix}\right]\,.
$$ Note that since $V^{kn}$ and $ V^{k-1,n}$ are projections of the same variety $\overline{\mathcal A_\pi}$ and $V^{kn}$ has one more coordinate than $ V^{k-1,n}$, there are only two possibilities for their dimensions:  $\dim V^{kn}=\dim V^{k-1,n}$ or $\dim V^{kn}=\dim V^{k-1,n}+1$.

 Now we begin the proof. By induction on $n$. For $n=1$ the statement is obviously true. Let us consider an $n\times n$ rank-control matrix $R(\pi)$ for some
involution $\pi\in S_n$. Its upper-left $(n-1)\times (n-1)$
submatrix is the rank-control matrix $R(\pi_{n-1})$ for the
involution $\pi_{n-1}\in S_{n-1}$ which corresponds to the upper-left
$(n-1)\times (n-1)$ submatrix of the matrix
$\pi$. By the induction hypothesis,
$$
\dim\,\overline{\mathcal A_{\pi_{n-1}}}=\frac{(n-1)^2-(n-1)}{2}-\mathfrak A\left(\pi_{n-1}\right)\,.
$$
 Now we add the $n$-th column to the partial involution matrix (the $n$-th does not provide any new information because we deal with the anti-symmetric matrices) and consider the $n$-th column of $R(\pi)$. (We also add the $n$-th row but since our matrices are anti-symmetric it suffices to understand only what happens to the dimension when we add the $n$-th column.) We added $n-1$ new coordinates (since there must be zero at the position $(n,n)$) to the variety $\overline{\mathcal A_{\pi_{n-1}}}$ and we have to show that
$$
\dim\,\overline{\mathcal A_{\pi}}=\dim\,\overline{\mathcal A_{\pi_{n-1}}}+n-1-\#\left\{(i,n)\,|\,1\leqslant i\leqslant n-1\,\,\textrm{and}\,\, r_{in}=r_{i-1,n-1}\right\}\,,\eqno(*)
$$
 i.e. not all the $n-1$ coordinates that we added make the dimension greater but only those of them  for which there is an inequality in the corresponding place of the certain diagonal of the rank-control matrix. The equality $(*)$ implies the statement of our theorem since $\frac{(n-1)^2-(n-1)}{2}+n-1=\frac{n^2-n}{2}$ and
$$
\mathfrak A\left(\pi\right)=\mathfrak A\left(\pi_{n-1}\right)+\#\left\{(i,n)\,|\,1\leqslant i\leqslant n-1\quad\textrm{and}\quad r_{in}=r_{i-1,n-1}\right\}\,.
$$
 Obviously, if $r_{1,n}=0$, then $a_{1,n}=0$ for any
$A=\left(a_{ij}\right)^n_{i,j=1}\in\overline{\mathcal A_\pi}$, this
itself is a polynomial equation which makes the dimension lower by
1, while if $r_{1,n}=1$ it means that the rank of the first row is
maximal and therefore, no equation. (In other words,  the dimension of the variety $V^{1n}$ which corresponds to $\left[\begin{matrix}
0 &a_{12} &\cdots &a_{1,n-1} &a_{1,n}\\
-a_{12} &0 &\cdots &a_{2,n-1} &\square\\
\cdots &\cdots &\cdots &\cdots &\square\\
-a_{1,n-1} &-a_{2.n-1} &\cdots &0  &\square\\
-a_{1,n} &\square &\square &\square &0
\end{matrix}\right]$ is greater by one than the of the variety $V^{0n}$ which corresponds to $\left[\begin{matrix}
0 &a_{12} &\cdots &a_{1,n-1}\\
-a_{12} &0 &\cdots &a_{2,n-1} \\
\cdots &\cdots &\cdots &\cdots \\
-a_{1,n-2} &-a_{2,n-2} &\cdots &a_{n-2,n-1} \\
-a_{1,n-1} &-a_{2,n-1} &\cdots &0
\end{matrix}\right]$ when $r_{1,n}=1$ and they have equal dimensions when $r_{1,n}=0$.) Now move down along the $n$-th
column of $R(\pi)$. Again by induction, this time the induction is
on the number of row $k$, assume that for each $1\leqslant i\leqslant k-1$ the dimensions of $V^{in}$ and $V^{i-1,n}$ ar equal iff
$r_{i-1,n-1}=r_{i,n}$ and  $\dim V^{in}=\dim V^{i-1,n}+1$ iff $r_{i-1,n-1}<r_{i,n}$.
First, let $r_{k-1,n-1}=r_{k,n}=c$. Consider a matrix
$A=\left(a_{ij}\right)_{i,j=1}^n\in\overline{\mathcal A_\pi}$ and
consider its upper-left $(k-1)\times(n-1)$ submatrix
$\left[\begin{matrix}
0 &a_{12} &\cdots &\cdots &a_{1,n-1}\\
-a_{12} &0 &\cdots &\cdots &a_{2,n-1}\\
\cdots &\cdots &\cdots &\cdots &\cdots \\
-a_{1,k-1} &-a_{2,k-1} &\cdots &\cdots &a_{k-1,n-1}
\end{matrix}\right]$. Using the notation introduced in
Proposition~\ref{ulr}, we denote this submatrix as $A_{k-1,n-1}$.  If $c=0$, then $rank A_{kn}=0$, so $A_{kn}$ is a zero matrix and $\dim V^{in}=\dim V^{i-1,n}=0$. Let $c\neq0$.
Since $ rank\left(A_{k-1,n-1}\right)=c$, we can take $c$ linearly independent columns
$\left[\begin{matrix}a_{1,j_1}\\a_{2,j_1}\\ \cdots\\
a_{k-1,j_{1}}\end{matrix}\right]$ , ... , $\left[\begin{matrix}a_{1,j_c}\\a_{2,j_c}\\ \cdots\\
a_{k-1,j_c}\end{matrix}\right]$ which span its column space. Now take only linearly independent rows of the $(k-1)\times c$ matrix $\left[\begin{matrix}
a_{1,j_1} &\cdots &a_{1,j_c} \\
a_{2,j_1} &\cdots &a_{2,j_c} \\
\cdots &\cdots &\cdots\\
a_{k-1,j_1} &\cdots &a_{k-1,j_c}\end{matrix}\right]$ to get a nonsingular $c\times c$ matrix $\left[\begin{matrix}
a_{i_1,j_1} &\cdots &a_{i_1,j_c} \\
a_{i_2,j_1} &\cdots &a_{i_2,j_c} \\
\cdots &\cdots &\cdots\\
a_{i_c,j_1} &\cdots &a_{i_c,j_c}\end{matrix}\right]$. The
equality $r_{k-1,n-1}=r_{k,n}=c\leqslant k-1$ implies that any $(c+1)\times(c+1)$ minor of the matrix $A_{kn}$ is zero, in particular
$\det\left[\begin{matrix}
a_{i_1,j_1} &\cdots &a_{i_1,j_c} &a_{i_1,n}\\
a_{i_2,j_1} &\cdots &a_{i_2,j_c} &a_{i_2,n}\\
\cdots &\cdots &\cdots &\cdots\\
a_{i_c,j_1} &\cdots &a_{i_c,j_c} &a_{i_c,n}\\
a_{k,j_1} &\cdots &a_{k,j_c} &a_{k,n}
\end{matrix}\right]=0$, which is a polynomial equation. This equation is algebraically independent of the similar equations obtained for $1\leqslant i\leqslant k-1$ because it involves the "new" variable -- the entry $a_{k,n}$. It indeed
 involves the entry $a_{k,n}$ since $\det\left[\begin{matrix}
a_{i_1,j_1} &\cdots &a_{i_1,j_c} \\
a_{i_2,j_1} &\cdots &a_{i_2,j_c} \\
\cdots &\cdots &\cdots \\
a_{i_c,j_1} &\cdots &a_{i_c,j_c}
\end{matrix}\right]\neq 0$. This equation means that the variable $a_{k,n}$ is not independent of the coordinates of the variety $V^{k-1,n}$, and therefore $\dim V^{k-1,n}=\dim V^{kn}$.

Now  let $r_{k-1,n-1}<r_{k,n}=c$, and we have to show that in this
case the variable $a_{nk}$ is independent of the coordinates of $V^{k-1,n}$, in other words, we have to show that there is no new equation. Consider the fragment
$\left[\begin{matrix}
r_{k-1,n-1} &r_{k-1,n}\\
r_{k-1,n} &r_{k,n}
\end{matrix}\right]$. There are four possible cases:
\begin{align*}
\left[\begin{matrix}
r_{k-1,n-1} &r_{k-1,n}\\
r_{k-1,n} &r_{k,n}
\end{matrix}\right]=&
\left[\begin{matrix}
c-1 &c-1\\
c-1  &c
\end{matrix}\right]\quad\textrm{or}\quad
\left[\begin{matrix}
c-2 &c-1\\
c-1  &c
\end{matrix}\right]
\quad\textrm{or}\\
& \left[\begin{matrix}
c-1 &c\\
c-1  &c
\end{matrix}\right]\quad\textrm{or}\quad
\left[\begin{matrix}
c-1 &c-1\\
c  &c
\end{matrix}\right]\,\,.
\end{align*}
The equality $r_{k,n}=c$ implies that each $(c+1)\times(c+1)$ minor  of $A_{kn}$ is equal to zero, but we shall see that each such equation is not new, i.e. it is implied by  the equality $r_{k,n-1}=c-1$ or by the equality $r_{k-1,n}=c-1$. In the first three of above four cases we decompose the $(c+1)\times(c+1)$ determinant $\det\left[\begin{matrix}
\cdots &\cdots   \\
\cdots  &a_{k,n}
\end{matrix}\right]$ using the last column. Since in all these cases $r_{k,n-1}=c-1$, each $c\times c$ minor of this decomposition (i.e. each $c\times c$ minor of $A_{k,n-1}$) is zero and therefore, this determinant is zero. In the fourth case we get the same if we decompose the determinant using its last row instead of the last column: since $r_{k-1,n}=c-1$, all the $c\times c$ minors of this decomposition (i.e. all $c\times c$ minor of $A_{k-1,n}$) are zeros and thus, our $(c+1)\times(c+1)$ determinant equals to zero. So there is no algebraic dependence between $a_{kn}$ and the coordinates of $V^{k-1,n}$. Therefore, $\dim V^{kn}=\dim V^{k-1,n}+1$.
  The proof is completed.
\end{proof}

Notice that the number $\frac{n^2-n}{2}$ is the dimension of the vector space of $n\times n$ anti-symmetric matrices.

To illustrate Theorem~\ref{PosetRankFunction} let us consider the second (from the top) level in the example given below. The diagonals (with added zeros at the beginning) of the upper triangle of the rank-control matrix $R_1=\left[\begin{matrix}
0 &1 &1 &1\\
1 &2 &2 &2\\
1 &2 &2 &2\\
1 &2 &2 &2 \end{matrix}\right]$ are: $(0\,\,1\,\,2\,\,2)$, $(0\,\,1\,\,2)$ and $(0\,\,1)$. We see that there is only one equality in the first of them (we have twice 2) and therefore $\mathfrak A(R_1)=1$. The diagonals (with added zeros at the beginning) of the upper triangle of the rank-control matrix $R_2=\left[\begin{matrix}
0 &0 &1 &1\\
0 &0 &1 &2\\
1 &1 &2 &3\\
1 &2 &3 &4
 \end{matrix}\right]$ are: $(0\,\,0\,\,1\,\,3)$, $(0\,\,1\,\,2)$ and $(0\,\,1)$. We see that there is only one equality in the first of them (we have twice 0) and therefore $\mathfrak A(R_2)=1$. The fact that $\mathfrak A(R_1)=\mathfrak A(R_2)=1$ means that the matrices $R_1$ and $R_2$ are both at the first from the top level in the poset.

 \subsection{An example.} Here we present this poset for $n=4$. In the first diagram we have the monomial matrices which are the representatives of the orbits and corresponding involutions of $S_4$. In the second diagram there are rank-control matrices.
$$
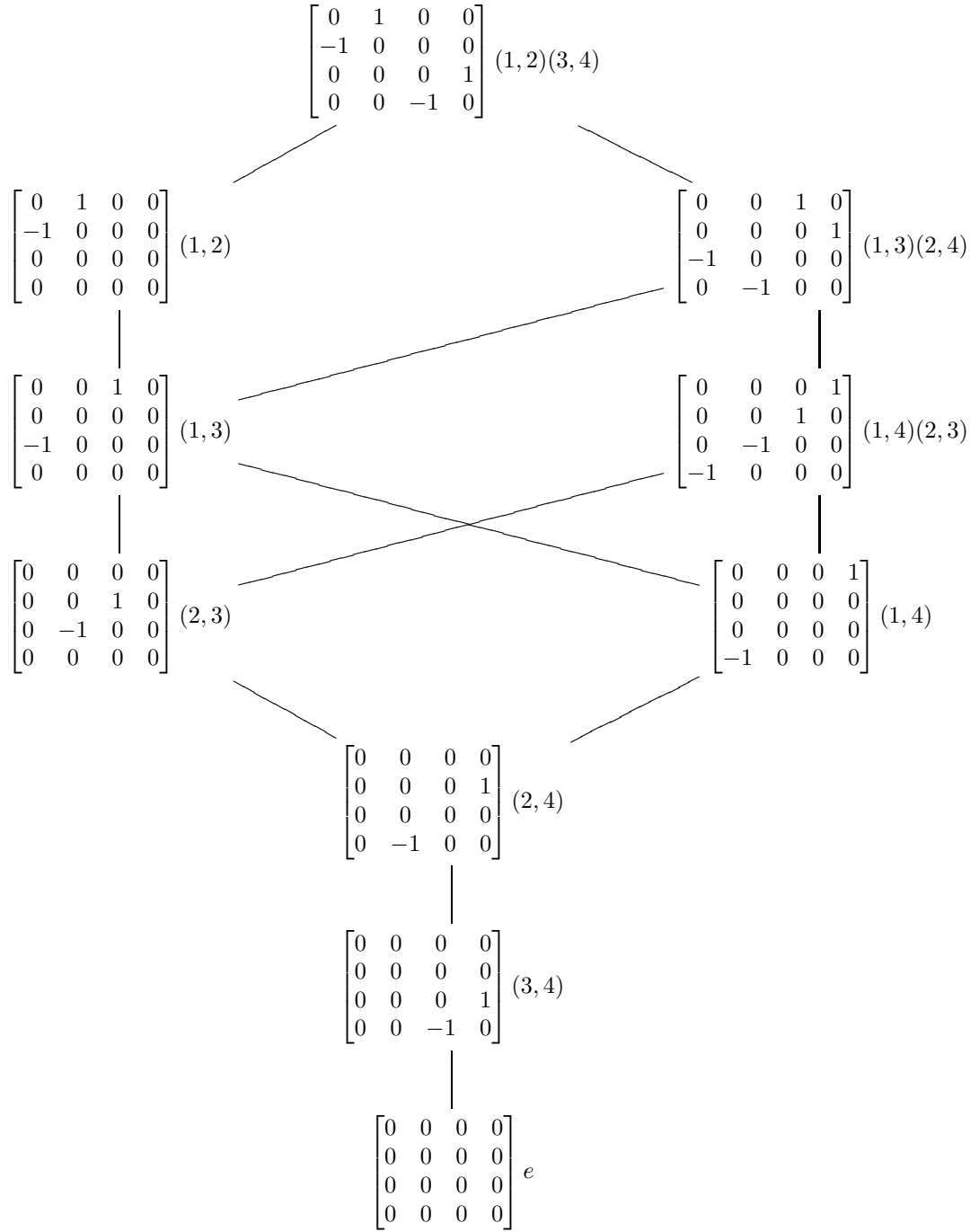
\xymatrix{& {\left[\begin{matrix}
0 &1 &0 &0\\
-1 &0 &0 &0\\
0 &0 &0 &1\\
0 &0 &-1 &0 \end{matrix}\right](1,2)(3,4)}& \\
{\left[\begin{matrix}
 0 &1 &0 &0\\
-1 &0 &0 &0\\
0 &0 &0 &0\\
0 &0 &0 &0 \end{matrix}\right](1,2)}\ar@{-}[ur] &  &{\left[\begin{matrix}
0 &0 &1 &0\\
0 &0 &0 &1\\
-1 &0 &0 &0\\
0 &-1 &0 &0
\end{matrix}\right](1,3)(2,4)} \ar@{-}[ul]\\
{\left[\begin{matrix}
0 &0 &1 &0\\
0 &0 &0 &0\\
-1 &0 &0 &0\\
0 &0 &0 &0
\end{matrix}\right](1,3)}\ar@{-}[u]\ar@{-}[urr] &  &{\left[\begin{matrix}
0 &0 &0 &1\\
0 &0 &1 &0\\
0 &-1 &0 &0\\
-1 &0 &0 &0   \end{matrix}\right](1,4)(2,3)} \ar@{-}[u]\\
{\left[\begin{matrix}
0 &0 &0 &0\\
0 &0 &1 &0\\
0 &-1 &0 &0\\
0 &0 &0 &0  \end{matrix}\right](2,3)}\ar@{-}[u]\ar@{-}[urr] &  &{\left[\begin{matrix}
0 &0 &0 &1\\
0 &0 &0 &0\\
0 &0 &0 &0\\
-1 &0 &0 &0 \end{matrix}\right](1,4)} \ar@{-}[u]\ar@{-}[ull]\\
& {\left[\begin{matrix}
0 &0 &0 &0\\
0 &0 &0 &1\\
0 &0 &0 &0\\
0 &-1 &0 &0\end{matrix}\right](2,4)}\ar@{-}[ul]\ar@{-}[ur]\\
 & {\left[\begin{matrix}
0 &0 &0 &0\\
0 &0 &0 &0\\
0 &0 &0 &1\\
0 &0 &-1 &0\end{matrix}\right](3,4)}\ar@{-}[u]\\
 & {\left[\begin{matrix}
0 &0 &0 &0\\
0 &0 &0 &0\\
0 &0 &0 &0\\
0 &0 &0 &0\end{matrix}\right]e}\ar@{-}[u]}
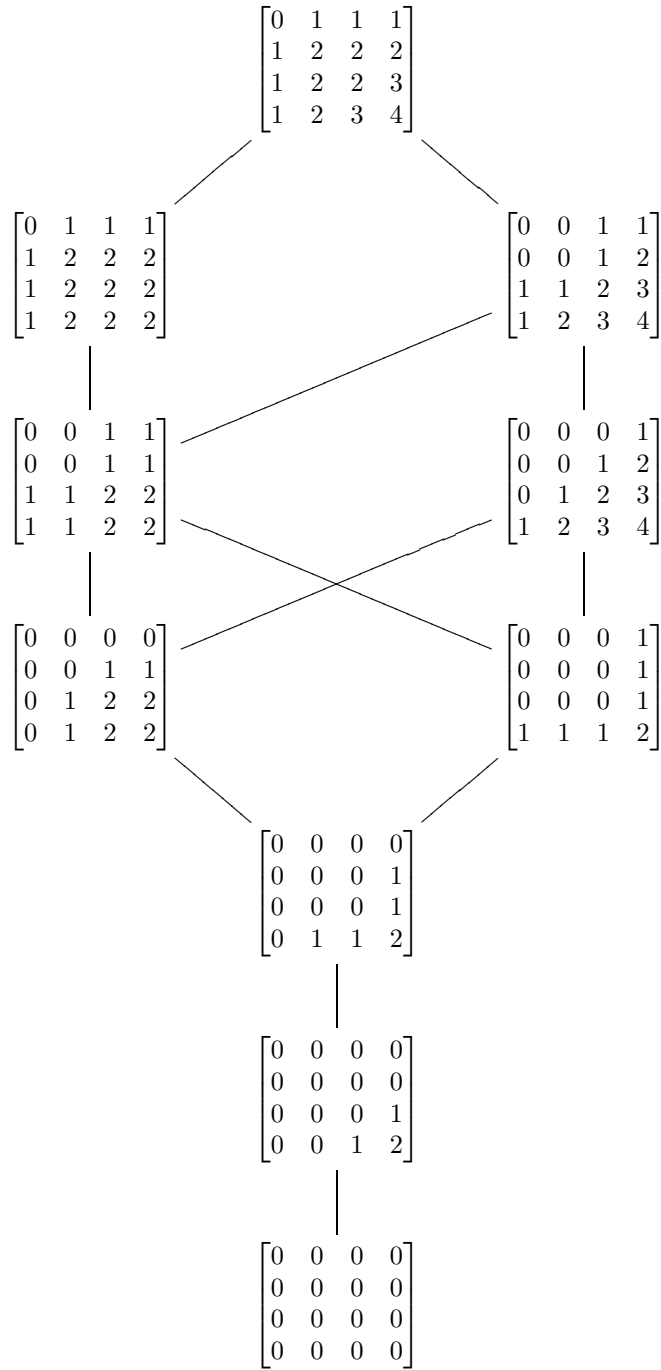
$$
$$\xymatrix{& {\left[\begin{matrix}
0 &1 &1 &1\\
1 &2 &2 &2\\
1 &2 &2 &3\\
1 &2 &3 &4 \end{matrix}\right]}& \\
{\left[\begin{matrix}
0 &1 &1 &1\\
1 &2 &2 &2\\
1 &2 &2 &2\\
1 &2 &2 &2 \end{matrix}\right]}\ar@{-}[ur] &  &{\left[\begin{matrix}
0 &0 &1 &1\\
0 &0 &1 &2\\
1 &1 &2 &3\\
1 &2 &3 &4
 \end{matrix}\right]} \ar@{-}[ul]\\
{\left[\begin{matrix}
0 &0 &1 &1\\
0 &0 &1 &1\\
1 &1 &2 &2\\
1 &1 &2 &2 \end{matrix}\right]}\ar@{-}[u]\ar@{-}[urr] &  &{\left[\begin{matrix}
0 &0 &0 &1\\
0 &0 &1 &2\\
0 &1 &2 &3\\
1 &2 &3 &4
 \end{matrix}\right]} \ar@{-}[u]\\
{\left[\begin{matrix}
0 &0 &0 &0\\
0 &0 &1 &1\\
0 &1 &2 &2\\
0 &1 &2 &2  \end{matrix}\right]}\ar@{-}[u]\ar@{-}[urr] &  &{\left[\begin{matrix}
0 &0 &0 &1\\
0 &0 &0 &1\\
0 &0 &0 &1\\
1 &1 &1 &2 \end{matrix}\right]} \ar@{-}[u]\ar@{-}[ull]\\
& {\left[\begin{matrix}
0 &0 &0 &0\\
0 &0 &0 &1\\
0 &0 &0 &1\\
0 &1 &1 &2\end{matrix}\right]}\ar@{-}[ul]\ar@{-}[ur]\\
 & {\left[\begin{matrix}
0 &0 &0 &0\\
0 &0 &0 &0\\
0 &0 &0 &1\\
0 &0 &1 &2\end{matrix}\right]}\ar@{-}[u]\\
 & {\left[\begin{matrix}
0 &0 &0 &0\\
0 &0 &0 &0\\
0 &0 &0 &0\\
0 &0 &0 &0\end{matrix}\right]}\ar@{-}[u]}
$$
\section{Another formula for the rank function.}
In this section, as before, we also don't distinguish between an involution $\pi\in S_n$ and the monomial anti-symmetric matrix (with minuses in the lower triangle) associated to $\pi$ by the bijection presented in Proposition~\ref{asorball}.
\bde Let $\pi\in S_n$ be an involution. It is always possible to write it as product of disjoint transpositions
$$
\pi=\left(i_1,j_1\right)\left(i_2,j_2\right)\cdots\left(i_k,j_k\right)
$$
in such a way that for all $1\leqslant t\leqslant k$, $i_t<j_t$ and $i_1<i_2<\cdots<i_k$. Let us call it "the canonic form".

Denote by $\mathfrak I(\pi)$ the {\bf number of inversions} in the word $i_1j_1i_2j_2\cdots i_kj_k$.
\ede

\bpr\label{SecFm} Let $\pi\in S_n$ be an involution. Then
$$
\mathfrak A(\pi)=\mathfrak I(\pi)+\sum_{a\,:\,\pi(a)=a}(n-a)\,.
$$
\epr
\begin{proof} By induction on $n$. Denote by $\pi_{n+1}$ some involution of $S_{n+1}$ and also the monomial anti-symmetric matrix which corresponds to this involution. Denote by $\pi_n$ the involution of $S_n$ which corresponds to the upper-left $n\times n$ block of the matrix $\pi_{n+1}$. For $\pi_n$ the statement is true by the induction hypothesis. Now we have to consider two cases:\\
{\bf Case1.} The number $n+1$ is a fixed point of the involution $\pi_{n+1}$ or in other words the $n$-th row and column of the matrix $\pi_{n+1}$ consist only of zeros. In this case we obviously have $\mathfrak I(\pi_{n+1})=\mathfrak I(\pi_n)$ and each fixed of $\pi_n$ is also a fixed point of $\pi_{n+1}$. It means that each fixed point of $\pi_n$ contributes an additional 1 (with comparison with $\mathfrak A(\pi_n$)) to $\mathfrak A(\pi_{n+1}$. So,
\begin{align*}
\mathfrak A(\pi_{n+1})&=\mathfrak A(\pi_{n})+\sum_{a\,:\,\pi_n(a)=a}1=\\
&=\mathfrak I(\pi_n)+\sum_{a\,:\,\pi_n(a)=a}(n-a)+\sum_{a\,:\,\pi_n(a)=a}1=\\
&=\mathfrak I(\pi_{n+1})+\sum_{a\,:\,\pi_{n+1}(a)=a}(n+1-a)\,.
\end{align*}
{\bf Case2.} The number $n+1$ is not a fixed point of the involution $\pi_{n+1}$. It means that the canonic form of $\pi_{n+1}$ is
$$
\pi_{n+1}=\left(i_1,j_1\right)\cdots(i,n+1)\cdots\left(i_k,j_k\right)
$$
for some $i$ or in other words the matrix $\pi_{n+1}$ has 1 at the position $(i,n+1)$ (and,of course, it also has $-1$ at the position $(n+1,i)$). We must show that
$$
\mathfrak A(\pi_{n+1})=\mathfrak I(\pi_{n+1})+\sum_{a\,:\,\pi_{n+1}(a)=a}(n+1-a)\,.
$$
By obvious observation
$$
\mathfrak A(\pi_{n+1})=\mathfrak A(\pi_{n})+\#\left\{a\,:\,\pi_n(a)=a\,\&\,a<i\right\}\,.
$$
By induction hypothesis $\mathfrak A(\pi_{n})=\mathfrak I(\pi_n)+\sum_{a\,:\,\pi_n(a)=a}(n-a)$ and so
\begin{equation*}
\mathfrak A(\pi_{n+1})=\mathfrak I(\pi_n)+\sum_{a\,:\,\pi_n(a)=a}(n-a)+\#\left\{a\,:\,\pi_n(a)=a\,\&\,a<i\right\}\,.\eqno{(*)}
\end{equation*}
Now, by direct calculation we have
\begin{align*}
\mathfrak I(\pi_{n+1})&+\sum_{a\,:\,\pi_{n+1}(a)=a}(n+1-a)=\mathfrak I(\pi_{n})+\#\left\{\left(i_t,j_t\right)\,:\,i_t<i\,\&\,j_t>i\right\}+\\
&+
2\cdot\#\left\{\left(i_t,j_t\right)\,:\,i_t>i\right\}+\sum_{a\,:\,\pi_{n+1}(a)=a}(n-a)
+\#\left\{a\,:\,\pi_{n+1}(a)=a\right\}=\\
&=\mathfrak I(\pi_{n})+\#\left\{\left(i_t,j_t\right)\,:\,i_t<i\,\&\,j_t>i\right\}+
2\cdot\#\left\{\left(i_t,j_t\right)\,:\,i_t>i\right\}+\\
&+\sum_{a\,:\,\pi_n(a)=a}(n-a)-(n-i)+\#\left\{a\,:\,\pi_{n+1}(a)=a\right\}\qquad\qquad\,\,\,\,\,\,\,\qquad(**)
\end{align*}
Comparing $(*)$ and $(**)$ we see that suffices to show that
\begin{align*}
\#\left\{a\,:\,\pi_n(a)=a\,\&\,a<i\right\}&=
\#\left\{\left(i_t,j_t\right)\,:\,i_t<i\,\&\,j_t>i\right\}+
2\cdot\#\left\{\left(i_t,j_t\right)\,:\,i_t>i\right\}+\\
&+\#\left\{a\,:\,\pi_{n+1}(a)=a\right\}-(n-i)
\end{align*}
which is obviously equivalent to the equality
\begin{align*}
n-i =& \#\left\{\left(i_t,j_t\right)\,:\,i_t<i\,\&\,j_t>i\right\}+
2\cdot\#\left\{\left(i_t,j_t\right)\,:\,i_t>i\right\}+\\
&+
\#\left\{a\,:\,\pi_n(a)=a\,\&\,a>i\right\}\,,
\end{align*}
which is indeed true since the number $$\#\left\{\left(i_t,j_t\right)\,:\,i_t<i\,\&\,j_t>i\right\}+
2\cdot\#\left\{\left(i_t,j_t\right)\,:\,i_t>i\right\}$$ is the number of all numbers greater than $i$ which are not fixed by $\pi_{n+1}$ while the number $\#\left\{a\,:\,\pi_n(a)=a\,\&\,a>i\right\}$ is the number of all numbers greater than $i$ which are fixed by $\pi_{n+1}$. So the sum of these two numbers is indeed $n-i$.
\end{proof}
\section{The final remark}
It follows from Observation~\ref{asorbinv} and the results of~\cite{ABC} that the poset of congruence B-orbits of invertible anti-symmetric matrices which is the poset of involutions not having fixed points in the symmetric group is a graded subposet of the (reversed) Bruhat poset of involutions of the symmetric group studied by F.\ Incitti in~\cite{I}, while the whole poset of involutions of the symmetric group with the ordering given by the containment of closures of congruence B-orbits of anti-symmetric matrices is not the Bruhat poset. The same is true not only for the set of involutions without fixed points but for any set of involutions with prescribed support (or in other words fixed set of fixed points).

Notice that the subposet of involutions with prescribed support is not an interval in the Bruhat poset of involutions of the symmetric group while it is an interval in the poset of involutions of the symmetric group introduced in this paper.\\
\noindent
{\bf Acknowledgements.} I would like to express the special gratitude to Prof. Lex Renner for the very useful information about the Bruhat poset and answering the numerous questions. This work was done when I was a postdoc at the Department of Mathematics at the University of Geneva and I would like to thank the department and especially Prof. Tatiana Smirnova-Nagnibeda
for the hospitality and discussion. I am grateful to the Swiss National Science Foundation for the financial support. I am grateful to Prof. Ron M. Adin, Dr. Eli Bagno, Dr. Anna Melnikov and Prof. Yuval Roichman  for the helpful discussions.

\end{document}